\documentclass[12pt]{amsart}
\usepackage{amssymb}
\newtheorem{theorem}{Theorem}[section]
\newtheorem{lemma}[theorem]{Lemma}
\newtheorem{corollary}[theorem]{Corollary}
\newtheorem{proposition}[theorem]{Proposition}
\theoremstyle{definition}
\newtheorem{remark}[theorem]{Remark}
\newtheorem{definition}[theorem]{Definition}
\newtheorem{example}[theorem]{Example}
\numberwithin{equation}{section}
\newcommand{\A}{\mbox{${\mathcal A}$}}
\newcommand{\M}{\mbox{$\mathbb{M}$}}
\newcommand{\C}{\mbox{$\mathbb{C}$}}
\newcommand{\K}{\mbox{$\mathbb{K}$}}
\newcommand{\Hi}{\mbox{${\mathcal H}$}}
\newcommand{\B}{\mbox{$\mathbb{B}$}}
\newcommand{\N}{\mbox{$\mathbb{N}$}}

\newcommand{\Ch}{\mbox{${\mathcal C}$}}
\newcommand{\D}{\mbox{${\mathcal D}$}}
\newcommand{\E}{\mbox{${\mathcal E}$}}
\newcommand{\Z}{\mbox{$\mathbb{Z}$}}
\newcommand{\R}{\mbox{$\mathbb{R}$}}
\begin{document}
\title[Open projections and Murray-von Neumann equivalence]{Open projections and Murray-von Neumann equivalence}
\author{Masayoshi Kaneda}
\author{Thomas Schick}
\address{Mathematisches Institut\\Georg-August-Universit\"{a}t G\"{o}ttingen\\Bunsenstra{\ss}e 3\\D-37073 G\"{o}ttingen\\Deutschland}
\email{masayoshi.kaneda@mathematik.uni-goettingen.de}
\address{Mathematisches Institut\\Georg-August-Universit\"{a}t G\"{o}ttingen\\Bunsenstra{\ss}e 3\\D-37073 G\"{o}ttingen\\Deutschland}
\email{thomas.schick@math.uni-goettingen.de}
\date{\today}
\begin{abstract}We characterize the $C^\star$-algebras for which openness of projections in their second duals is preserved under Murray-von Neumann equivalence. They are precisely the extensions of the annihilator $C^\star$-algebras by the commutative $C^\star$-algebras.

We also show that the annihilator $C^\star$-algebras are precisely the $C^\star$-algebras for which all projections in their second duals are open.
\end{abstract}
\maketitle
\section{Introduction}\label{section:intro}
The notion of open (respectively, closed, compact) projections was introduced by Akemann as a noncommutative generalization of open (respectively, closed, compact) subsets of a topological space, and a remarkable theory of ``noncommutative topology'' was developed including noncommutative versions of the Stone-Weierstrass theorem, the Urysohn lemma, etc. in a series of his pioneer works \cite{Akemann1969}, \cite{Akemann1970}, \cite{Akemann1971}. He called the collection of open projections the \emph{hull-kernel structure} (\emph{HKS} for short) which can be considered as a ``noncommutative topology.''

Let us recall the definitions. Throughout the paper, our $C^\star$-algebras are not necessarily unital. Let $\A$ be a $C^\star$-algebra and $\A^{**}$ be its second dual. We recall that the second dual $\A^{**}$ of a $C^\star$-algebra $\A$ is canonically a $W^\star$-algebra, and we regard $\A$ as a $C^\star$-subalgebra of $\A^{**}$ in the canonical way. As a $W^\star$-algebra, $\A^{**}$ has the canonical type decomposition $\A^{**}=A^{**}_{\operatorname{I}}\oplus\A^{**}_{\operatorname{II}}\oplus\A^{**}_{\operatorname{III}}$. A projection $p\in\A^{**}$ is \emph{open} if there exists an increasing net $(a_\alpha)$ of positive elements in $\A$ such that $a_\alpha\nearrow p$ in the weak$^*$ topology on $\A^{**}$, or, equivalently, if $(p\A^{**}p\cap\A)^{\perp\perp}=p\A^{**}p$. Throughout the paper a \emph{projection} means an orthogonal projection, that is, $p=p^\star=p^2$. A projection $q\in\A^{**}$ is \emph{closed} if $1_{\scriptsize\A^{**}}-q$ is open. A projection $r\in\A^{**}$ is \emph{compact} if it is closed and there exists an $a\in\A$ with $a\ge0$ such that $ra=r$. Projections $p,q\in\A^{**}$ are said to be \emph{Murray-von Neumann equivalent} (or, simply, \emph{equivalent}) and we write $p\sim q$ if there exists a partial isometry $v\in\A^{**}$ such that $p=vv^\star$ and $q=v^\star v$. We write $p\precsim q$ if $p$ is Murray-von Neumann equivalent to a subprojection of $q$.

In preparing \cite{Kaneda-qmideal}, the first author encountered the problem ``Is every projection which is Murray-von Neumann equivalent to an open projection in $\A^{**}$ open?'' The answer turns out to be negative as the following example shows.
\begin{example}\emph{There exist a $C^{\star}$-algebra $\A$ and Murray-von Neumann equivalent projections $p,q\in\A^{**}$ such that $p$ is open but $q$ is not.}
\end{example}
\begin{proof}Let $\A:=C([0,1])\otimes\M_2$, where $\M_2$ is the $C^\star$-algebra of $2\times2$ matrices with entries in $\C$. Recall that $\ell^\infty([0,1])$ is canonically contained in $C([0,1])^{**}$. Let $F:=\{0\}$ and $F^c:=(0,1]$, then the characteristic functions $\chi_F$ and $\chi_{F^c}$ are in $\ell^\infty([0,1])$. Define a partial isometry $v:=\chi_F\otimes e_{11}+\chi_{F^c}\otimes e_{12}\in\A^{**}$, where $e_{ij}\in\M_2$ is the canonical matrix unit, that is, it has $1$ in its $(i,j)$-entry and $0$ elsewhere. Then the projection $p:=vv^{\star}=1\otimes e_{11}$ is open, whereas $q:=v^{\star}v=\chi_F\otimes e_{11}+\chi_{F^c}\otimes e_{22}$ is not.
\end{proof}Let us define openness of projections to be \emph{stable} under Murray-von Neumann equivalence in $\A^{**}$ if for two Murray-von Neumann equivalent projections in $\A^{**}$ openness of one implies openness of the other.

The above example motivates us to ask the question ``For which $C^\star$-algebras is openness of projections in their second duals stable under Murray-von Neumann equivalence?'' Two obvious kinds of such $C^\star$-algebras are the commutative $C^\star$-algebras and the $C^\star$-algebras for which every projection in their second duals is open.

Although Murray-von Neumann equivalence is an important notion, especially in the classification of von Neumann algebras and $K$-theory for $C^\star$-algebras, it is not so strong as to preserve many structures determined by projections. For instance, the hereditary $C^\star$-subalgebras corresponding to two Murray-von Neumann equivalent open projections need not be $^\star$-isomorphic \cite[Theorem~9]{Lin1990}. Furthermore, in recent years various kinds of equivalence relations between open projections have been studied in classifying $C^\star$-algebras and in connection with Cuntz semigroups (\cite{PZ2000}, \cite{Lin2010}, \cite{ORT2011}, \cite{NW2015-1}; also see \cite{Cuntz1978}). None of them, however, discuss stability of openness under an equivalence relation, so it will be important to consider the above question. The main result (Theorem~\ref{th:main}) of this paper characterizes the $C^\star$-algebras for which openness of projections in their second duals is stable under Murray-von Neumann equivalence, and it turns out that the aforementioned two kinds of $C^\star$-algebras are ``essentially'' the only ones.
\begin{theorem}\label{th:main}Let $\A$ be a $C^{\star}$-algebra. Then openness of projections in $\A^{**}$ is stable under Murray-von Neumann equivalence if and only if there exists an ideal $J$ in $\A$ such that $J$ is an annihilator $C^\star$-algebra and the quotient $C^\star$-algebra $\A/J$ is commutative.
\end{theorem}In conjunction with the main result we also show that the annihilator $C^\star$-algebras are precisely those $C^\star$-algebras for which all projections in their second duals are open (Theorem~\ref{th:allopen}). This result does also follow from \cite{ABS2012}, although it is not explicitly stated there. We provide our own direct proof without relying on results in \cite{ABS2012}, which we think is of independent interest.

We use the following definition of an annihilator $C^\star$-algebra convenient to us.
\begin{definition}An \emph{annihilator $C^\star$-algebra} is a $C^\star$-algebra isomorphic to a $c_0$-direct sum $\oplus_{i\in I}^{c_0}\K(\Hi_i)$ for some index set $I$, where $\K(\Hi_i)$ is the $C^\star$-algebra of compact operators on a Hilbert space $\Hi_i$. (See the end of this section for the definition of a $c_0$-direct sum.)
\end{definition}
\begin{remark}Historically, Kaplansky called such algebras \emph{dual}, but we prefer to call them annihilator $C^\star$-algebras following \cite{ABS2012} to avoid confusion with Banach space dual. There are numerous characterizations of the annihilator $C^\star$-algebras (e.g.~\cite[Section 2]{Kaplansky1951}, \cite[4.7.20]{Dixmier1977}, \cite[Theorem~5.5]{ABS2012}), and Theorem~\ref{th:allopen} adds another to the list of characterizations. Annihilator $C^\star$-algebras are also called \emph{compact $C^\star$-algebras} in the literature, but again we do not use this terminology to avoid confusion with a single summand $\K(\Hi)$.
\end{remark}We exclusively use the symbol $\cong$ to indicate that the two $C^\star$-algebras connected by it are $^\star$-isomorphic. As a well-known fact, any $^\star$-isomorphism between $W^\star$-algebras (von Neumann algebras) is also a homeomorphism with respect to the weak$^\star$ ($\sigma$-weak operator) topologies concerned. Let $J$ be an ideal in a $C^\star$-algebra $\A$. In our paper, an \emph{ideal} always means a norm-closed two-sided ideal. Recall that an exact sequence $0\to J\to\A\to\A/J\to 0$ induces an exact sequence $0\to J^{**}\to\A^{**}\to(\A/J)^{**}\to0$ (as taking the second duals is an exact functor), so that $\A^{**}/J^{**}\cong(\A/J)^{**}$ canonically. Since $J^{\perp\perp}\,(\cong J^{**})$ is a weak$^*$-closed ideal in $\A^{**}$, there exists a central projection $z\in\A^{**}$ such that $J^{\perp\perp}=z\A^{**}$, so that the last exact sequence splits, hence $\A^{**}\cong J^{**}\oplus\A^{**}/J^{**}$.

If $\{p_\alpha\}$ is a set of projections in a $W^\star$-algebra $M$, then the \emph{meet} $\bigwedge_\alpha p_\alpha$ is defined to be the largest projection majorized by every $p_\alpha$, and the \emph{join} $\bigvee_\alpha p_\alpha$ is the smallest projection majorizing every $p_\alpha$. If $N$ is a $W^\star$-subalgebra of $M$ and $\{p_\alpha\}\subseteq N$, then the meet and the join defined in $N$ are the same as the meet and the join defined in $M$.

The symbol $\otimes_{\min}$ stands for the minimum (also called spatial, or injective) $C^\star$-tensor product, and $\overline{\otimes}$ stands for the $W^\star$-tensor product. For a locally compact Hausdorff space $\Omega$, $C_0(\Omega):=\{f\in C(\Omega)\,|\,\forall\varepsilon>0,\{\omega\in\Omega\,|\,|f(\omega)|\ge\varepsilon\}\text{ is compact}\}$. If $\Omega$ is discrete, we prefer to write $c_0(\Omega)$ instead. Finally, for normed spaces $X_i$ indexed by a set $I$, their \emph{$c_0$-direct sum} is defined to be
\begin{equation*}\bigoplus_{i\in I}^{c_0}X_i:=\left\{x\in\prod_{i\in I}X_i\,\middle|\,\forall\varepsilon>0,\,\{i\in I\,|\,\|x(i)\|\ge\varepsilon\}\text{ is a finite set}\right\}.
\end{equation*}
\section{Open projections and Murray-von Neumann equivalence}\label{section:main}
This section is devoted to prove the main result, Theorem~\ref{th:main}, of this paper. We need the following proposition as well as several lemmas.
\begin{proposition}\label{pr:continuous part}Let $\A$ be a $C^{\star}$-algebra. Then the continuous part (i.e., the non-type~I part) of $\A^{**}$ cannot contain a nonzero open projection.
\end{proposition}
\begin{proof}Suppose the contrary and let $p$ be a nonzero open projection contained in the continuous part of $\A^{**}$, and put $\A_p:=p\A^{**}p\cap\A$. Then $p\A^{**}p$, which is $^\star$-isomorphic to $(\A_p)^{**}$, must have a nonzero type~I summand (e.g. \cite[the first paragraph]{Chu1991}), a contradiction.
\end{proof}The following lemma directly follows from the definition of open projections.
\begin{lemma}\label{lm:join is open}Suppose that $\A$ is a $C^{\star}$-algebra and that $\{p_\alpha\}$ is a set of open projections in $\A^{**}$. Then $\bigvee_\alpha{p_\alpha}$ is open.
\end{lemma}The following lemma is a special case of \cite[Theorem~II.7]{Akemann1969} in terms of open projections. In \cite{Akemann1969}, however, $C^{\star}$-algebras are assumed to be unital, so a justification is in order.
\begin{lemma}\label{lm:p meets q}Suppose that $\A$ is a $C^{\star}$-algebra and that $p_1,p_2,p\in\A^{**}$ are mutually orthogonal projections such that $q_1:=p_1+p$ and $q_2:=p_2+p$ are open. Then $p$ is also open.
\end{lemma}
\begin{proof}Let $\A^1$ be the unitization of $\A$, and we denote by $1$ the identity of $\A^1$. By \cite[Theorem~2.4]{BHN2008}, a projection $q\in\A^{**}$ is open in $\A^{**}$ if and only if it is open in $(\A^1)^{**}$. As a consequence of \cite[Theorem~II.7]{Akemann1969}, the assumptions imply that $p$ is open in $(\A^1)^{**}$ and hence in $\A^{**}$.
\end{proof}
\begin{remark}Open projections in $\A^{**}$ remain open in $(\A^1)^{**}$, but closed projections in $\A^{**}$ need not be closed in $(\A^1)^{**}$. Therefore it is easier to generalize \cite[Proposition~II.5 and Theorem~II.7]{Akemann1969} to the nonunital case in terms of open projections as in the above lemma.
\end{remark}
\begin{lemma}\label{lem:pass_to_z}Let $\A$ be a $C^\star$-algebra such that openness of projections in $\A^{**}$ is stable under Murray-von Neumann equivalence. If every subprojection of a projection $p\in\A^{**}$ is open, then every subprojection of its central support $z(p)$ is also open.
\end{lemma}
\begin{proof}Let $q$ be a subprojection of $z(p)$. Set
\begin{equation*}Q:=\{q'\in\A^{**}\,|\,\text{$q'$ is a projection such that $q'\le q$ and $q'\precsim p$}\}.
\end{equation*}Then $q-\bigvee Q$ is a subprojection of $z(p)$. If it was non-zero, it would by \cite[Lemma~V.1.7]{Takesaki2002} contain a subprojection $q''$ equivalent to a subprojection of $p$. Then $q''\in Q$, a contradiction. Hence $q=\bigvee Q$, and it is open by Lemma~\ref{lm:join is open} since each element in $Q$ is open by stability.
\end{proof}The following is the key lemma.
\begin{lemma}\label{lm:key}Let $\A$ be a $C^\star$-algebra such that openness of projections in $\A^{**}$ is stable under Murray-von Neumann equivalence. Then every non-abelian central projection in $\A^{**}$ majorizes a nonzero central projection all of whose subprojections are open.
\end{lemma}
\begin{proof}Represent $\A$ universally on a Hilbert space $\Hi$, and identify its SOT-closure with $\A^{**}$. Let $z\in\A^{**}$ be a non-abelian central projection. Then it follows from \cite[Lemma~V.1.7]{Takesaki2002} that there exist nonzero projections $p_0,q_0\in\A^{**}$ such that $p_0\le z$, $q_0\le z$, $p_0\perp q_0$, and $p_0\sim q_0$. Let $v\in\A^{**}$ be a partial isometry implementing the equivalence, that is, $p_0=vv^\star$ and $q_0=v^\star v$. The idea now is to choose a self-adjoint $a\in\A$ which is SOT-close to $p_0$. This gives rise to the open spectral projections corresponding to $\sigma(a)\cap(1/2,\|a\|]$ and $\sigma(a)\cap[0,1/2)$. These spectral projections are the projections to which we apply the assumption that openness is preserved under Murray-von Neumann equivalence. More concretely, pick $\eta\in\Hi$ with $\|\eta\|=1$ such that $q_0\eta=\eta$, and put $\xi:=v\eta$. Then $p_0\xi=\xi$ and $\|\xi\|=1$. By Kaplansky's density theorem we choose a self-adjoint $a\in\A$ with $0\le a\le1$ such that $\|(p_0-a)\xi\|<1/8$ (hence $\|a\|\ge\|a\xi\|>7/8$) and $\|a\eta\|<1/8$. Let $e_1$ and $e_2$ be the spectral projections of $a$ in $\A^{**}$ corresponding to $\sigma(a)\cap(1/2,\|a\|]$ and $\sigma(a)\cap[0,1/2)$, respectively. Then $e_1$ and $e_2$ are open projections, and $e_1^\perp=\overline{e}_2$ and $e_2^\perp=\overline{e}_1$, where $\overline{e}_1$ is the closure of $e_1$, i.e., the meet of the closed projections majorizing $e_1$. Since $a\ge\frac{1}{2}\overline{e}_1$, we have that
\begin{equation*}\frac{1}{8}>\|a\eta\|\ge\langle a\eta,\eta\rangle\ge\frac{1}{2}\langle\overline{e}_1\eta,\eta\rangle=\frac{1}{2}\|\overline{e}_1\eta\|^2
\end{equation*}and hence
\begin{equation}\label{eq:eta}\|\overline{e}_1\eta\|<\frac{1}{2}.
\end{equation}Since $a\le\|a\|e_1+\frac{1}{2}\overline{e}_2$, we also have that
\begin{align*}\frac{3}{4}&<\left(\frac{7}{8}\right)^2<\|a\xi\|^2\le\|a^{\frac{1}{2}}\xi\|^2=\langle a\xi,\xi\rangle\le\|a\|\langle e_1\xi,\xi\rangle+\frac{1}{2}\langle\overline{e}_2\xi,\xi\rangle\\&=\|a\|\|e_1\xi\|^2+\frac{1}{2}\|\overline{e}_2\xi\|^2\le\|e_1\xi\|^2+\frac{1}{2}
\end{align*}and hence
\begin{equation}\label{eq:xi}\|e_1\xi\|>\frac{1}{2}.
\end{equation}Now we have that
\begin{equation*}\|e_1ve_2\eta\|=\|e_1v(\eta-\overline{e}_1\eta)\|=\|e_1\xi-e_1v\overline{e}_1\eta\|\ge\|e_1\xi\|-\|\overline{e}_1\eta\|>0,
\end{equation*}where we used Inequalities~(\ref{eq:xi})~and~(\ref{eq:eta}) in the last inequality. Let $p_1$ and $p_2$ be, respectively, the left and right support projections of $e_1ve_2$ in $\A^{**}$. Then $p_1\le ze_1$, $p_2\le ze_2$, $p_1\perp p_2$, and $p_1\sim p_2$. Let $w\in\A^{**}$ be a partial isometry implementing the equivalence $p_1\sim p_2$, that is $p_1=ww^\star$ and $p_2=w^\star w$. If $p$ is any subprojection of $p_1$, then $p\sim w^\star pw\le p_2\le e_2$, and $e_1\sim e_1-p+w^\star pw$, hence by stability $e_1-p+w^\star pw$ is open since $e_1$ is open. The projections $e_1-p$, $e_2-w^\star pw$, and $w^\star pw$ satisfy the hypotheses of Lemma~\ref{lm:p meets q}, hence $w^\star pw$ is open, and so is $p$ by stability. Hence every subprojection of $p_1$ is open. By Lemma \ref{lem:pass_to_z}, its central support projection $z(p_1)$ is the object we had to construct.
\end{proof}We are now in a position to prove Theorem~\ref{th:main}.
\begin{proof}[\textbf{\emph{Proof of Theorem~\ref{th:main}}}] The ``if'' direction easily follows from $\A^{**}\cong J^{**}\oplus(\A/J)^{**}$ since $(\A/J)^{**}$ is commutative and all projections in $J^{**}$ are open by Theorem~\ref{th:allopen}.

To show the converse, suppose that $\A$ is a $C^{\star}$-algebra such that openness of projections in $\A^{**}$ is stable under Murray-von Neumann equivalence. It follows from Lemma~\ref{lm:key} and Proposition~\ref{pr:continuous part} that $\A^{**}$ is of type~I. Let $z_{\operatorname{I}_1}\in\A^{**}$ be the central projection onto the type~I$_1$ part (i.e., the abelian type~I part). Then $z_{\operatorname{I}_1}^\perp\,(:=1_{\scriptsize\A^{**}}-z_{\operatorname{I}_1})$ does not majorize any abelian central projection. Put
\begin{align*}Z:=\{z\in\A^{**}\,|\,& \text{$z$ is an open central projection such that} \\ & \text{$z\le z_{\operatorname{I}_1}^\perp$ and every subprojection of $z$ is open}\}.
\end{align*}Then $\bigvee Z=z_{\operatorname{I}_1}^\perp$ since otherwise by Lemma~\ref{lm:key} $z_{\operatorname{I}_1}^\perp-\bigvee Z$, which is a non-abelian central projection, majorizes a nonzero open central projection all of whose subprojections are open, a contradiction. Let $p\in\A^{**}$ be any subprojection of $z_{\operatorname{I}_1}^\perp$. Then $p=\bigvee\{zp\,|\,z\in Z\}$ is open by the definition of $Z$ and Lemma~\ref{lm:join is open}. Put $J_0:=z_{\operatorname{I}_1}^\perp\A^{**}\cap\A$, which is an ideal in $\A$ and also an annihilator $C^\star$-algebra by Theorem~\ref{th:allopen}. Furthermore, $\A/J_0$ is commutative since its second dual $(\A/J_0)^{**}\cong z_{\operatorname{I}_1}\A^{**}$ is commutative.
\end{proof}
\begin{remark}
\begin{enumerate}
  \item In Theorem \ref{th:main} the ideal $J$ is not unique in general. This is due to the arbitrariness of including commutative annihilator $C^\star$-algebras in $J$. There exist, however, the largest one and the smallest one. The $J_0$ obtained in the proof above is the smallest one (which does not contain any commutative direct summand), and the largest one includes all the atomic (= minimal nonzero) projections in $\A$.
  \item A nonzero commutative direct summand of $\A^{**}$ can occur even if $\A$ has no commutative direct summand, as the following example shows.
\end{enumerate}
\end{remark}
\begin{example}Set $\Omega:=\{0\}\cup\{1/n\,|\,n\in\Z^+\}\subset\R$. Then $\Omega$ is a compact Hausdorff space. Set $\mathring{\Omega}:=\Omega\setminus\{0\}$, which is discrete. Let$$\A:=\{f\in C(\Omega,\M_2)\,|\,f(0)\text{ is supported on the $(1,1)$-entry}\}.$$Then $\A$ does not have a commutative direct summand, and $J$ is uniquely determined to be $J\cong C_0(\mathring{\Omega},\M_2)\cong\bigoplus_{\omega\in\mathring{\Omega}}^{c_0}\M_2$, while $\A^{**}$ has a commutative direct summand: $\A^{**}\cong\C\oplus\bigoplus_{\omega\in\mathring{\Omega}}^{\ell^\infty}\M_2$.
\end{example}
\appendix
\section*{Appendix. Noncommutative discreteness and annihilator $C^\star$-algebras}\label{section:discreteness}
\renewcommand{\thesection}{A}
\setcounter{theorem}{0}
\renewcommand{\thetheorem}{\thesection.\arabic{theorem}}
Since open projections are a replacement of open sets, one can imagine that the noncommutative notion corresponding to the discrete topology, in which every singleton is open, will be that every projection is open. The following proposition suggests that this is the correct idea.
\begin{proposition}\label{pr:allopen commutative}Let $\Omega$ be a locally compact Hausdorff space. Then every projection in $C_0(\Omega)^{**}$ is open if and only if $\Omega$ is discrete. In this case $f\in C_0(\Omega)$ is supported on an at-most-countable subset of $\Omega$, therefore $C_0(\Omega)$ is an ideal in $C_0(\Omega)^{**}$, and hence $C_0(\Omega)^{**}$ is the multiplier algebra of $C_0(\Omega)$. In particular, if $\Omega$ is compact, then $C_0(\Omega)\,(=C(\Omega))$ is $\ell_n^{\infty}$ for some finite number $n\in\N$.
\end{proposition}
\begin{proof}The ``if'' direction is clear since $C_0(\Omega)^{**}=\ell^{\infty}(\Omega)$ if $\Omega$ is discrete. For the converse,  if $\omega\in\Omega$ then the characteristic function $\chi_\omega$ is in $\ell^\infty(\Omega)$ which is canonically contained in $C_0(\Omega)^{**}$. By assumption $\chi_\omega$ is an open projection in $C_0(\Omega)^{**}$ which implies that $\{\omega\}$ is an open set.
\end{proof}We need two lemmas before proceeding to the general $C^\star$-algebra case.
\begin{lemma}\label{lm:compact}Let $\K(\Hi)$ be the $C^{\star}$-algebra of compact operators on a (possibly nonseparable) Hilbert space $\Hi$. Then every projection in $\K(\Hi)^{**}$ is open.
\end{lemma}
\begin{proof}Every rank-$1$ projection is in $\K(\Hi)$, and every projection in $\K(\Hi)^{**}\,(\cong\B(\Hi))$ is a join of rank-$1$ projections, thus it is open by Lemma~\ref{lm:join is open}.
\end{proof}
\begin{lemma}\label{lm:minproj}Let $\A$ be a $C^{\star}$-algebra. If a projection $0\ne p\in\A^{**}$ is minimal and open, then $p\in\A$.
\end{lemma}
\begin{proof}Let $p\ne0$ be a minimal projection in $\A^{**}$. Then by \cite[Proposition~5.1]{BN2012} it is compact. Hence by \cite[Theorem~2.2~(iii)$\Rightarrow$(ii)]{BN2012} it is closed in $(\A^1)^{**}$, where $\A^1$ is the unitization of $\A$. If $p$ is open in $\A^{**}$, it remains open in $(\A^1)^{**}$, so it is clopen in $(\A^1)^{**}$. Then by \cite[Proposition~II.18]{Akemann1969} $p\in\A^1$, hence $p\in\A^{\perp\perp}\cap\A^1=\A$.
\end{proof}
\begin{theorem}\label{th:allopen}Let $\A$ be a $C^{\star}$-algebra. Then every projection in $\A^{**}$ is open if and only if $\A$ is an annihilator $C^\star$-algebra.
\end{theorem}
\begin{proof}Suppose that $\A$ is an annihilator algebra, i.e., $\A\cong\bigoplus_{k\in I}^{c_0}\K(\Hi_k)$. Then $\A^{**}\cong(\bigoplus_{k\in I}^{c_0}\K(\Hi_k))^{**}\cong\bigoplus_{k\in I}^{\ell^\infty}\K(\Hi_k)^{**}$. Since all projections in $\K(\Hi_k)^{**}$ are open by Lemma~\ref{lm:compact}, all projections in $\bigoplus_{k\in I}^{\ell^\infty}\K(\Hi_k)^{**}$ are open by Lemma~\ref{lm:join is open}.

To show the converse, suppose that every projection in $\A^{**}$ is open. By Proposition~\ref{pr:continuous part} $\A^{**}$ must be of type~I. By the structure theory of type~I $W^{\star}$-algebras (e.g. \cite[III.1.5.12]{Blackadar2006}) there exist a set $I$ of cardinals and mutually orthogonal nonzero central projections $\{z_n\,|\,n\in I\}$ in $\A^{**}$ such that $\sum_{n\in I}z_n=1$ (the identity of $\A^{**}$) and $\A^{**}z_n$ is of type~I$_n$, and furthermore, each $\A^{**}z_n$ is $^{\star}$-isomorphic to $\B(\Hi_n)\overline{\otimes}\D_n$, where $\Hi_n$ is an $n$-dimensional Hilbert space and $\D_n$ is a commutative $W^{\star}$-algebra.

Put $\A_n:=\A^{**}z_n\cap\A$, then $\A^{**}z_n=\A_n^{\perp\perp}\,(\cong\A_n^{**})$ since $z_n$ is open in $\A^{**}$ by assumption. Noting that by \cite[Theorem~2.4]{BHN2008} a projection in $\A_n^{**}$ is open in $\A_n^{**}$ if and only if it is open in $\A^{**}$, we shall focus our argument for the moment on a ``summand'' $\A_n^{**}\cong\B(\Hi_n)\overline{\otimes}\D_n$ for an arbitrarily fixed $n\in I$. Henceforth we consider as $\A_n^{**}=\B(\Hi_n)\overline{\otimes}\D_n$ instead of $\A_n^{**}\cong\B(\Hi_n)\overline{\otimes}\D_n$ for notational simplicity.

Let $p$ be any rank-$1$ projection in $\B(\Hi_n)$ and put
\begin{equation*}\Ch_n:=(p\otimes1_{\scriptsize\D_n})\A_n^{**}(p\otimes1_{\scriptsize\D_n})\cap\A_n.
\end{equation*}Then $\Ch_n^{**}\cong\Ch_n^{\perp\perp}=(p\otimes1_{\scriptsize\D_n})\A_n^{**}(p\otimes1_{\scriptsize\D_n})$ since $p\otimes1_{\scriptsize\D_n}$ is open in $\A_n^{**}$ by assumption. But $(p\otimes1_{\scriptsize\D_n})\A_n^{**}(p\otimes1_{\scriptsize\D_n})=p\B(\Hi_n)p\otimes\D_n\cong\D_n$ since $p\B(\Hi_n)p=\C p$. Therefore $\D_n$ is $^\star$-isomorphic to the second dual of the commutative $C^\star$-algebra $\Ch_n$ which we shall $^\star$-isomorphically identify with $C_0(\Omega_n)$, where $\Omega_n$ is a locally compact Hausdorff space. Let $q$ be any projection in $C_0(\Omega_n)^{**}$ (considered as $=\D_n$). Then by assumption $p\otimes q$ is open in $\A_n^{**}$, hence by \cite[Theorem~2.4]{BHN2008} $p\otimes q$ is open in $\Ch_n^{**}$ as well. Since $p\otimes q$ is open for every projection $q\in C_0(\Omega_n)^{**}$, $\Omega_n$ is discrete by Proposition~\ref{pr:allopen commutative}. Thus we shall henceforth write $c_0(\Omega_n)$ instead of $C_0(\Omega_n)$.

Now let $\chi_\omega$ be the characteristic function of a singleton $\{\omega\}\subseteq\Omega_n$. Then $\chi_\omega$ is a projection in $c_0(\Omega_n)$. Put
\begin{equation*}\E_\omega:=(1_{\scriptsize\Hi_n}\otimes\chi_\omega)\A_n^{**}(1_{\scriptsize\Hi_n}\otimes\chi_\omega)\cap\A_n.
\end{equation*}Then $\E_\omega^{**}\cong\E_\omega^{\perp\perp}=(1_{\scriptsize\Hi_n}\otimes\chi_\omega)\A_n^{**}(1_{\scriptsize\Hi_n}\otimes\chi_\omega)=\B(\Hi_n)\otimes\C\chi_\omega$ since $1_{\scriptsize\Hi_n}\otimes\chi_\omega$ is open in $\A_n^{**}$ by assumption. For any rank-$1$ projection $p\in\B(\Hi_n)$, $p\otimes\chi_\omega$ is a minimal projection in $\E_\omega^{**}$, therefore $p\otimes\chi_\omega\in\E_\omega$ by Lemma~\ref{lm:minproj}. Thus $\K(\Hi_n)\otimes\C\chi_\omega\subseteq\E_\omega\subseteq\A_n$ since $\E_\omega$ is a $C^\star$-algebra. Since $c_0(\Omega_n)=\bigoplus_{\omega\in\Omega_n}^{c_0}\C\chi_\omega$, we also get that $\K(\Hi_n)\otimes_{\min}c_0(\Omega_n)\subseteq\A_n$, hence
\begin{equation*}J:=\bigoplus^{c_0}_{n\in I}(\K(\Hi_n)\otimes_{\min}c_0(\Omega_n))\subseteq\A.
\end{equation*}Now $J$ is an ideal in $\A$, hence there exists a central projection $z\in\A^{**}$ such that $J^{\perp\perp}=z\A^{**}$. Suppose that $z\ne1$. Then there exists $n\in I$ such that $(1-z)z_n\ne0$. Since $(1-z)z_n\ne0$ is a central projection in $\A_n^{**}$, it must be of the form $1_{\scriptsize\Hi_n}\otimes\chi$ for some nonzero characteristic function $\chi\in c_0(\Omega_n)$. Thus there exists an $\omega\in\Omega_n$ such that $\E_\omega\subseteq(1-z)\A_n^{**}$, hence $(1-z)\A_n^{**}\cap J^{\perp\perp}\ne\{0\}$, a contradiction. Therefore $z=1$, so that $J^{\perp\perp}=\A^{**}$, hence
\begin{equation*}\A=J\cong\bigoplus^{c_0}_{n\in I}\bigoplus^{c_0}_{\omega\in\Omega_n}\K(\Hi_n),
\end{equation*}as we had to prove.
\end{proof}
\begin{corollary}Let $\A$ be a unital $C^{\star}$-algebra. Then every projection in $\A^{**}$ is open if and only if $\A$ is finite-dimensional.
\end{corollary}
\begin{proof}Looking at its structure, an annihilator $C^\star$-algebra is unital if and only if it is finite dimensional.  On the other hand, by the known structure of finite dimensional $C^\star$-algebras they are all annihilator $C^\star$-algebras.
\end{proof}
\section*{Acknowledgments}The first author wishes to express his gratitude to the second author and Ralf Meyer for their invitation to the Georg-August-Universit\"{a}t G\"{o}ttingen in summer 2018, which initiated this work. The first author thanks Mathematisches Institut of the Georg-August-Universit\"{a}t G\"{o}ttingen for the financial support during his second visit to G\"{o}ttingen in spring 2019. The first author also thanks Heidi and Joseph for their mental support.

\end{document}